\documentclass[12pt]{article}
\usepackage{latexsym,amsfonts,amsmath,amssymb}
\newtheorem{theorem}{Theorem}
\newtheorem{corollary}[theorem]{Corollary}
\newtheorem{sublemma}{Lemma}[theorem]
\newtheorem{lemma}[theorem]{Lemma}
\newtheorem{question}[theorem]{Question}
\newtheorem{observation}[theorem]{Observation}
\newtheorem{claim}[theorem]{Claim}
\newtheorem{conjecture}[theorem]{Conjecture}
\newtheorem{definition}[theorem]{Definition}
\newtheorem{remark}[theorem]{Remark}
\newtheorem{example}[theorem]{Example}
\def\Theorem #1. #2\par{\setbox1=\hbox{#1}\ifdim\wd1=0pt
   \begin{theorem}#2\end{theorem}\else
   \newtheorem{#1}[theorem]{#1}\begin{#1}\label{#1}#2\end{#1}\fi}
\def\Corollary #1. #2\par{\setbox1=\hbox{#1}\ifdim\wd1=0pt
   \begin{corollary}#2\end{corollary}\else
   \newtheorem{#1}[theorem]{#1}\begin{#1}\label{#1}#2\end{#1}\fi}
\def\Lemma #1. #2\par{\setbox1=\hbox{#1}\ifdim\wd1=0pt
   \begin{lemma}#2\end{lemma}\else
   \newtheorem{#1}[theorem]{#1}\begin{#1}\label{#1}#2\end{#1}\fi}
\def\SubLemma #1. #2\par{\setbox1=\hbox{#1}\ifdim\wd1=0pt
   \begin{sublemma}#2\end{sublemma}\else
   \newtheorem{#1}{#1}[theorem]\begin{#1}\label{#1}#2\end{#1}\fi}
\def\Question #1. #2\par{\setbox1=\hbox{#1}\ifdim\wd1=0pt
   \begin{question}#2\end{question}\else
   \newtheorem{#1}[theorem]{#1}\begin{#1}\label{#1}#2\end{#1}\fi}
\def\Observation #1. #2\par{\setbox1=\hbox{#1}\ifdim\wd1=0pt
   \begin{observation}#2\end{observation}\else
   \newtheorem{#1}[theorem]{#1}\begin{#1}\label{#1}#2\end{#1}\fi}
\def\Claim #1. #2\par{\setbox1=\hbox{#1}\ifdim\wd1=0pt
   \begin{claim}#2\end{claim}\else
   \newtheorem{#1}[theorem]{#1}\begin{#1}\label{#1}#2\end{#1}\fi}
\def\Conjecture #1. #2\par{\setbox1=\hbox{#1}\ifdim\wd1=0pt
   \begin{conjecture}#2\end{conjecture}\else
   \newtheorem{#1}[theorem]{#1}\begin{#1}\label{#1}#2\end{#1}\fi}
\def\Definition #1. #2\par{\setbox1=\hbox{#1}\ifdim\wd1=0pt
   \begin{definition}#2\end{definition}\else
   \newtheorem{#1}[theorem]{#1}\begin{#1}\label{#1}#2\end{#1}\fi}
\def\Remark #1. #2\par{\setbox1=\hbox{#1}\ifdim\wd1=0pt
   \begin{remark}#2\end{remark}\else
   \newtheorem{#1}[theorem]{#1}\begin{#1}\label{#1}#2\end{#1}\fi}
\def\Example #1. #2\par{\setbox1=\hbox{#1}\ifdim\wd1=0pt
   \begin{example}#2\end{example}\else
   \newtheorem{#1}[theorem]{#1}\begin{#1}\label{#1}#2\end{#1}\fi}
\def\QuietTheorem #1. #2\par{\setbox1=\hbox{#1}\medskip\ifdim\wd1=0pt
   \proclaim{Theorem}{#2}\else\proclaim{#1}{#2}\fi}
\newcommand{\proclaim}[2]{\medskip\noindent{\bf #1 } {\sl#2}\par\medskip\noindent}
\def\Proclaim #1. #2\par{\proclaim{#1}{#2}}
\newenvironment{proof}{\noindent}{\kern2pt\QEDbox\par\bigskip}
\def\Proof#1: {\begin{proof}\setbox1=\hbox{#1}\ifdim\wd1=0pt{\bf Proof: }\else{\bf #1: }\fi}
\newcommand{\QED}{\end{proof}}

\def\Abstract #1\par{\medskip\begin{quotation}{\singlespaced\footnotesize{\noindent{\bf Abstract.~}#1}}\end{quotation}\medskip}
\def\Title #1\par{\title{#1}\maketitle}
\def\Author #1\par{\author{#1}}
\def\Section #1\par{\section{#1}}
\def\SubSection #1\par{\subsection{#1}}
\def\SubSubSection #1\par{\subsubsection{#1}}
\def\MidTitle #1\par{\bigskip\goodbreak\centerline{\small\bf #1}\bigskip}

\newcommand{\singlespaced}{\baselineskip=15pt}
\renewcommand{\P}{{\mathbb P}}
\newcommand{\Q}{{\mathbb Q}}
\newcommand{\R}{{\mathbb R}}

\newcommand{\Qdot}{\dot\Q}

\newcommand{\Gtail}{G_{\!\scriptscriptstyle\rm tail}}

\newcommand{\Ptail}{\P_{\!\scriptscriptstyle\rm tail}}

\newfont{\msam}{msam10 at 12pt}

\newcommand{\of}{\subseteq}

\newcommand{\elesub}{\prec}

\newcommand{\cof}{\mathop{\rm cof}}

\newcommand{\Th}{\mathop{\rm Th}}

\newcommand{\plus}{{+}}

\newcommand{\satisfies}{\models}

\newcommand{\possible}{\mathop{\raisebox{-1pt}{$\Diamond$}}}
\newcommand{\necessary}{\mathop{\raisebox{3pt}{\framebox[6pt]{}}}}
\newcommand{\cross}{\times}

\newcommand{\intersect}{\cap}

\newcommand{\card}[1]{{\left|#1\right|}}

\newcommand{\UnderTilde}[1]{{\setbox1=\hbox{$#1$}\baselineskip=0pt\vtop{\hbox{$#1$}\hbox to\wd1{\hfil$\sim$\hfil}}}{}}
\newcommand{\Undertilde}[1]{{\setbox1=\hbox{$#1$}\baselineskip=0pt\vtop{\hbox{$#1$}\hbox to\wd1{\hfil$\scriptstyle\sim$\hfil}}}{}}
\newcommand{\undertilde}[1]{{\setbox1=\hbox{$#1$}\baselineskip=0pt\vtop{\hbox{$#1$}\hbox to\wd1{\hfil$\scriptscriptstyle\sim$\hfil}}}{}}
\newcommand{\UnderdTilde}[1]{{\setbox1=\hbox{$#1$}\baselineskip=0pt\vtop{\hbox{$#1$}\hbox to\wd1{\hfil$\approx$\hfil}}}{}}
\newcommand{\Underdtilde}[1]{{\setbox1=\hbox{$#1$}\baselineskip=0pt\vtop{\hbox{$#1$}\hbox to\wd1{\hfil$\scriptstyle\approx$\hfil}}}{}}

\newcommand{\st}{\mid}

\newcommand{\iso}{\cong}
\def\<#1>{\langle\,#1\,\rangle}

\newcommand{\QEDbox}{\fbox{}}

\newcommand{\ORD}{\mathop{\hbox{\sc ord}}}
\newcommand{\ZFC}{\hbox{\sc zfc}}

\newcommand{\CH}{\hbox{\sc ch}}

\newcommand{\AD}{\hbox{\sc ad}}
\newcommand{\PD}{\hbox{\sc pd}}
\newcommand{\MA}{\hbox{\sc ma}}

\newcommand{\MP}{\hbox{\sc mp}}
\newcommand{\MPtilde}{\UnderTilde{\MP}}
\newcommand{\factordiagramup}[6]{$$\begin{array}{ccc}
#1&\raise3pt\vbox{\hbox to60pt{\hfill$\scriptstyle
#2$\hfill}\vskip-6pt\hbox{$\vector(4,0){60}$}}&#3\\ \vbox
to30pt{}&\raise22pt\vtop{\hbox{$\vector(4,-3){60}$}\vskip-22pt\hbox
to60pt{\hfill$\scriptstyle #4\qquad$\hfill}}
     &\ \ \lower22pt\hbox{$\vector(0,3){45}$}\ {\scriptstyle #5}\\
\vbox to15pt{}&&#6\\
\end{array}$$}
\newcommand{\factordiagram}[6]{$$\begin{array}{ccc}
#1&&\\ \ \ \raise22pt\hbox{$\vector(0,-3){45}$}\ {\scriptstyle #2}
&\raise22pt\hbox{$\vector(2,-1){90}$}\raise5pt\llap{$\scriptstyle#3$\qquad\quad}&\vbox
to25pt{}\\ #4&\raise3pt\vbox{\hbox to90pt{\hfill$\scriptstyle
#5$\hfill}\vskip-6pt\hbox{$\vector(4,0){90}$}}&#6\\
\end{array}$$}
\newcommand{\df}{\it} % use italic for definition terms

\begin{document}
\author{Joel David Hamkins\thanks{My research has been
supported by grants from the PSC-CUNY Research Foundation and the
NSF. I would like to thank James Cummings, Ali Enayat, Ilijas
Farah, Philip Welch and W. Hugh Woodin for helpful discussions
concerning various topics arising in this article, as well as
Paul Larson for introducing me to the topic and Christophe
Chalons for his
insightful email correspondence.}\\
\normalsize\sc The City University of New
York\thanks{Specifically, The College of Staten Island of CUNY
and The CUNY Graduate Center.}\\
\normalsize\sc Carnegie Mellon University\thanks{I am currently
on leave as Visiting Associate Professor at Carnegie Mellon
University, and I would like to thank the CMU Department of
Mathematical Sciences for their hospitality.}\\
{\footnotesize http://math.gc.cuny.edu/faculty/hamkins}}

\Title A simple maximality principle

\Abstract In this paper, following an idea of Christophe Chalons,
I propose a new kind of forcing axiom, the {\it Maximality
Principle}, which asserts that any sentence $\varphi$ holding in
some forcing extension $V^\P$ and all subsequent extensions
$V^{\P*\Qdot}$ holds already in $V$. It follows, in fact, that
such sentences must also hold in all forcing extensions of $V$.
In modal terms, therefore, the Maximality Principle is expressed
by the scheme
$(\possible\necessary\varphi)\implies\necessary\varphi$, and is
equivalent to the modal theory $S5$. In this article, I prove
that the Maximality Principle is relatively consistent with \ZFC.
A boldface version of the Maximality Principle, obtained by
allowing real parameters to appear in $\varphi$, is
equiconsistent with the scheme asserting that $V_\delta\elesub V$
for an inaccessible cardinal $\delta$, which in turn is
equiconsistent with the scheme asserting that $\ORD$ is Mahlo.
The strongest principle along these lines is
$\necessary\MPtilde$, which asserts that $\MPtilde$ holds in $V$
and all forcing extensions. From this, it follows that $0^\#$
exists, that $x^\#$ exists for every set $x$, that projective
truth is invariant by forcing, that Woodin cardinals are
consistent and much more. Many open questions remain.

\Section The Maximality Principle

Christophe Chalons has introduced a delightful new axiom,
asserting in plain language that anything that is forceable and
not subsequently unforceable is true. Specifically, in
\cite{Chalons1} he proposes the following principle:

\QuietTheorem The Chalons Maximality Principle. Suppose $\varphi$
is a sentence and for every inaccessible cardinal $\kappa$ there
is a forcing notion $\P\in V_\kappa$ such that
$V_\kappa^\P\satisfies\varphi$ and for any such $\P$ and any
further forcing $\Qdot$ it holds that
$V_\kappa^{\P*\Qdot}\satisfies\varphi$, then there are unboundedly
many inaccessible cardinals $\kappa$ with
$V_\kappa\satisfies\varphi$.

In this paper, I would like to present a more streamlined version
of Chalons' axiom, along with an equiconsistency analysis of the
principle in various forms. My intention with these streamlined
axioms, which I call the Maximality Principles, is to take the
essence of Chalons' idea to heart, and accord more closely with
the plain language slogan, {\it ``anything forceable and not
subsequently unforceable is true.''}

So let me begin. Define that a sentence $\varphi$ is {\df
possible} or {\df forceable} when it holds in some forcing
extension $V^\P$. The sentence $\varphi$ is {\df necessary} when
it holds in all forcing extensions $V^\P$ (by trivial forcing,
this includes $V$). Thus, $\varphi$ is {\df forceably necessary}
(or {\df possibly necessary}) when $\varphi$ holds in some
forcing extension $V^\P$ and all subsequent extensions
$V^{\P*\Qdot}$, i.e. when it is forceable that $\varphi$ is
necessary.

The forceably necessary sentences are therefore those which can be
turned `on' in a permanent way, so that no further forcing can
ever turn them off. Perhaps the simplest example of a forceably
necessary sentence is $V\not=L$, since this is easy to force, and
once it is forced, one cannot force its negation. More interesting
examples include the assertion ``{\it there is a real with
minimal constructibility degree over $L$}\,,'' which can be
forced via Sacks forcing---and once true, this assertion persists
to any forcing extension---and the assertion ``{\it there are a
proper class of inaccessible cardinals, if any,}'' which I
consider in Theorem \ref{IfAny} below. Conversely, the sentences
\CH\ or $\neg\CH$ are forceable but never necessary because they
can be forced true or false over any model of set theory. The
principle I have in mind, or at least the first approximation to
it, is the following:

\QuietTheorem Maximality Principle\footnote{Chalons
\cite{Chalons2} considers another principle approaching this, a
scheme that asserts that if a statement $\varphi$ has the feature
that for every inaccessible cardinal $\kappa$ there is a poset
$\P\in V_\kappa$ with $V_\kappa^\P\satisfies\varphi$ and for every
such $\P$ and any further forcing $\Qdot\in V_\kappa$ we have
$V_\kappa^{\P*\Qdot}\satisfies\varphi$, then $\varphi$ holds (in
$V$). Using this scheme as an intermediary, Chalons establishes
the consistency of the Chalons Maximality Principle by first
establishing the consistency of it, assuming the existence of a
measurable cardinal. The ideas of the proof of Theorem \ref{MP} in
this paper show that one can get away with less; indeed, to obtain
the consistency of Chalons' principle a Mahlo cardinal is more
than sufficient, and actually it suffices to have $V_\delta\elesub
V$ for an inaccessible cardinal $\delta$. Apart from this
question, though, I prefer to untie the core idea of the
expression {\it everything forceable and not unforceable is true}
from the question of whether the statements are true in or
forceable over every $V_\kappa$, and this seems to lead to the
more direct principle I have stated here.} (\MP). Any statement in
the language of set theory that is forceably necessary is true.

\noindent Here, I intend to assert $\MP$ as a scheme, the scheme
that asserts of every statement $\varphi$ in the language of set
theory that $$\hbox{\sl if $\varphi$ is forceably necessary, then
$\varphi$.}$$ I mean, of course, that $\varphi$ is true in $V$.
But it actually follows from this scheme that $\varphi$ must also
be true in every forcing extension of $V$, that is, that
$\varphi$ is necessary; this is proved below in Theorem
\ref{EquivalentForms}.

By Tarski's theorem on the non-definability of truth, there is in
general no first order way to express whether a statement is
forceable or true; and so there seems to be little hope of
expressing the Maximality Principle scheme in a single
first-order sentence of set theory. Thus, I remain content with
the Maximality Principle expressed as a scheme.

As I hope my terminology suggests, the Maximality Principle admits
a particularly natural expression in modal logic. Specifically, if
we take $\possible\varphi$ to mean that $\varphi$ is possible,
that is, that $\varphi$ is forceable, that it holds in some
forcing extension, and $\necessary\varphi$ to mean that $\varphi$
is necessary, that is, that $\varphi$ holds in all forcing
extensions, then we have a very natural Kripke frame
interpretation of modal logic in which the possible worlds are
simply the models of \ZFC\ and one world is accessible from
another if it is a forcing extension of that world.\footnote{A
{\df Kripke frame} allows one to provide semantics for a modal
theory. It consists of a collection of worlds and an
accessibility relation between the worlds. In such a frame, an
assertion is {\df possible} in a world if it is true in a world
accessible by that world; the assertion is {\df necessary} in a
world if it is true in all worlds accessible by that world. The
accessibility relation of this article, the one relating any
model of set theory with its forcing extensions, is both
reflexive and transitive. Thus, it is what is known as an $S4$
frame.}

In this Kripke frame, one can easily check many elementary modal
facts. The assertion
$\neg\possible\varphi\iff\necessary\neg\varphi$, for example,
expresses the trivial fact that $\varphi$ is not forceable if and
only if $\neg\varphi$ holds in all forcing extensions;
furthermore, its validity provides a duality between $\possible$
and $\necessary$. The assertion
$\necessary(\varphi\implies\psi)\implies
(\necessary\varphi\implies\necessary\psi)$ is similarly easy to
verify in this context. Further, since every model of set theory
is a (trivial) forcing extension of itself, we can conclude
$\varphi\implies\possible\varphi$, asserting that any true
statement is forceable (by trivial forcing), as well as the dual
version $\necessary\varphi\implies\varphi$, asserting that any
necessary statement is true. Using the fact that a forcing
extension of a forcing extension is a forcing extension of the
ground model, we can verify
$\possible\possible\varphi\implies\possible\varphi$, asserting
that any statement forceable over a forcing extension is
forceable over the ground model, and the dual version
$\necessary\varphi\implies\necessary\necessary\varphi$, asserting
that if a statement holds in all forcing extensions, then this
remains true in any forcing extension.

The axioms mentioned in the previous paragraph are exactly the
axioms constituting the modal theory known as $S4$ (see, for
example, \cite{Modal}). And since the modal operators $\possible$
and $\necessary$, as I have defined them in this article, are
expressible in the language of set theory, we can view these $S4$
axioms as simple theorems of \ZFC.

In this modal terminology, the assertion that $\varphi$ is
forceably necessary is exactly expressed by the assertion
$\possible\necessary\varphi$. The Maximality Principle $\MP$,
therefore, asserts of every $\varphi$ that
$$(\possible\necessary\varphi)\implies\varphi.$$ This way of
stating the axiom naturally leads one to the apparently stronger
principle, known as the Euclidean Axiom in modal logic:
$$(\possible\necessary\varphi)\implies\necessary\varphi.$$ As a
scheme, however, this apparently stronger principle is in fact no
stronger than the original principle, a fact I alluded to
earlier. In fact, the Maximality Principle can be stated in a
variety of equivalent forms:

\Theorem Equivalent Forms Theorem. The following schemes are
equivalent, where $\varphi$ ranges over all sentences in the
language of set theory: \label{EquivalentForms}
 \begin{enumerate}
 \item Every forceably
    necessary statement is true: $\possible\necessary\varphi\implies\varphi$.
 \item Every
    forceably necessary statement is necessary: $\possible\necessary\varphi\implies\necessary\varphi$.
 \item Every forceable statement is necessarily forceable (i.e. forceable over every forcing extension):
    $\possible\varphi\implies\necessary\possible\varphi$.
 \item Every true
    statement is necessarily forceable: $\varphi\implies\necessary\possible\varphi$.
 \item Every non-necessary statement is necessarily non-necessary:
 $\neg\necessary\varphi\implies\necessary\neg\necessary\varphi$.
 \end{enumerate}

\Proof: This is an elementary exercise in modal reasoning. Scheme
2 is actually a special case of scheme 1, obtained by replacing
$\varphi$ with $\necessary\varphi$ in 1 to obtain
$\possible\necessary\necessary\varphi\implies\necessary\varphi$.
This is just
$\possible\necessary\varphi\implies\necessary\varphi$ because
$\necessary\necessary\varphi$ is equivalent to
$\necessary\varphi$. Conversely, scheme 2 implies scheme 1
because $\necessary\varphi$ implies $\varphi$. Scheme 3 is the
contrapositive form of scheme 2 (with $\neg\varphi$ replacing
$\varphi$) and vice versa, because
$\neg\possible\necessary\neg\varphi$ is equivalent to
$\necessary\possible\varphi$ and $\neg\necessary\neg\varphi$ is
equivalent to $\possible\varphi$. Scheme 4, similarly, is the
contrapositive form of scheme 1. Finally, scheme 5 is obtained
from scheme 3 (applied to $\neg\varphi$) by pushing the negation
through to $\varphi$, and conversely by pulling it out again.\QED

Schemes 2, 3 and 5 can be strengthened to the full equivalences
$$\possible\necessary\varphi\iff\necessary\varphi,\qquad
\possible\varphi\iff\necessary\possible\varphi\quad\hbox{and}\quad
\neg\necessary\varphi\iff\necessary\neg\necessary\varphi$$
because the converse implication in each case is immediate.

Because the various schemes are equivalent, let me now freely
refer to any of them as the Maximality Principle \MP. In
particular, in form 2 above the principle asserts that any
possibly necessary statement is necessary, and thus, under the
Maximality Principle the collection of statements holding
necessarily is maximal.

Those readers familiar with modal logic will recognize Scheme 5
above, the axiom
$\neg\necessary\implies\necessary\neg\necessary\varphi$, as
precisely the axiom one adds to $S4$ to make the modal theory
$S5$ (see \cite{Modal}). Therefore, since we have already observed
that the weaker $S4$ axioms hold automatically, we can make the
following conclusion:

\Theorem. In the Kripke frame whose worlds are the models of
\ZFC, each of which accesses precisely its forcing extensions,
the Maximality Principle \MP\ in a world is equivalent to the
assertion of the $S5$ axioms of modal logic in that
world.\footnote{Let me emphasize here that the Maximality
Principle includes the assertions
$\possible\necessary\varphi\implies\necessary\varphi$ only for
{\it sentences} $\varphi$. In particular, it would be wrong to
include, as one might want in other contexts, the universal
generalizations of this axiom applied to formulas with free
variables. Indeed, doing so would result in an inconsistent
theory, as I prove in Observation \ref{Parameters}.}

One naturally wants to express the Maximality Principle in a
language with parameters, and so I will also consider the
following boldface version of the Maximality Principle, which
allows for arbitrary parameters from $\R$:

\QuietTheorem Maximality Principle, Boldface Version ($\MPtilde$).
Any statement in the language of set theory with arbitrary real
parameters that is forceably necessary is necessary. That is,
$(\possible\necessary\varphi)\implies\necessary\varphi$.

This boldface version of the principle also admits the equivalent
forms of Theorem \ref{EquivalentForms}. Subsuming the parameters
into a universal quantifier, one can alternatively think of the
boldface Maximality Principle as the scheme asserting all
sentences of the form
$$\forall x\in\R\,[\possible\necessary\varphi(x)\implies
\necessary\varphi(x)],$$
where $\varphi$ is a formula with one free variable $x$.

With real parameters, of course, one can in effect refer to any
hereditarily countable set as a parameter. And while one might
want to include even more parameters than this, it would be
inadvisable to do so in view of the following:

\Observation. If a set $p$ is not hereditarily countable, then
there is a statement $\varphi(p)$ in the language of set theory
with parameter $p$ which is forceably necessary, but not true.
Hence, the Maximality Principle asserted in the language of set
theory with arbitrary parameters is false.\label{Parameters}

\Proof: Let $\varphi(p)$ simply be the statement, ``{\it $p$ is
hereditarily countable}.'' Of course, this is forceably necessary,
since any set can be made countable by forcing, and once
countable, it can never be made uncountable again by further
forcing. But by the hypothesis on $p$, the statement $\varphi(p)$
is not true.\QED

Before introducing the next principle, let me observe the
following.

\Observation. The Maximality Principle \MP, if true, is necessary.
In modal terminology,
$\MP\implies\necessary\MP$.\label{MPimpliesBoxMP}

\Proof: If a sentence $\varphi$ is forceably necessary over a
forcing extension, then it was already forceably necessary over
the ground model, and if necessary in the ground model, it
remains necessary in any forcing extension. That is,
$\possible\varphi$ is downward absolute from a forcing extension,
and $\necessary\varphi$ is upward absolute to any forcing
extension. So $\MP\implies\necessary\MP$.\QED

Observation \ref{MPimpliesBoxMP} does not easily generalize to the
boldface Maximality Principle \MPtilde. The problem is that there
can be new real parameters in the forcing extension and thus
entirely new assertions $\varphi$ there, which cannot be
expressed in the ground model. This leads us to the strongest
principle I will consider in this article, the principle that
asserts:

\QuietTheorem Necessary Maximality Principle
($\necessary\MPtilde$). The principle $\MPtilde$ is necessary;
that is, it holds in $V$ and every forcing extension $V^\P$.

Let me caution the reader that one does not obtain the principle
$\necessary\MPtilde$ merely by prefacing every assertion in the
\MPtilde\ scheme with a $\necessary$\,, for this would only assert
the scheme in a forcing extension using parameters from the ground
model; such assertions hold for free by \MPtilde\ as in
Observation \ref{MPimpliesBoxMP}. Rather, $\necessary\MPtilde$
asserts that the full \MPtilde\ scheme holds in every forcing
extension $V^\P$, using the additional real parameters available
there.

Later, I will prove that the principle $\necessary\MPtilde$ is
much stronger by far than the other Maximality Principles. If it
holds, for example, then $0^\#$ exists and more, the universe is
closed under sharps and projective truth is invariant by set
forcing. An upper bound on the consistency strength of
$\necessary\MPtilde$ is not yet known; it may simply be false.

\Section Consistency of the Maximality Principle \MP

Let us now analyze the consistency of the simplest of the
Maximality Principles, namely, the lightface version \MP\ with no
parameters.

\Theorem. If\/ \ZFC\ is consistent, then so also is \ZFC\ plus the
Maximality Principle \MP.\label{MP}

\Proof: In order to highlight various aspects of the theory, I
will actually give two proofs of this theorem. The first proof is
elementary, but the second proof will generalize to $\MPtilde$.
We begin with a simple observation.

\SubLemma. For no sentence $\varphi$ are both $\varphi$ and
$\neg\varphi$ forceably necessary.

\Proof: Suppose that $\varphi$ is forced to be necessary by $\P$,
so that $\varphi$ holds in every forcing extension of $V^\P$, and
that $\neg\varphi$ is forced necessary by $\Q$, so that
$\neg\varphi$ holds in every forcing extension of $V^\Q$.
Consider now the forcing $\P\cross\Q$. Since this is isomorphic
to $\Q\cross\P$, the extension $V^{\P\cross\Q}$ is a forcing
extension both of $V^\P$ and of $V^\Q$. Thus, both $\varphi$ and
$\neg\varphi$ hold there, a contradiction.\QED

This argument generalizes to the following.\footnote{Modal
logicians will recognize that these lemmas rely on the {\it
directedness} of the accessibility relation on worlds.
Specifically, the crucial fact at play here is that any two
forcing extensions can be combined into a third, which is a
forcing extension of each of them.}

\SubLemma. Over any model of set theory, the collection of
forceably necessary statements forms a consistent
theory.\label{Consistent}

\Proof: Consider any finite collection
$\varphi_1,\varphi_2,\ldots,\varphi_n$ of such sentences, each
forceably necessary over a fixed model $M$. Suppose that
$\varphi_i$ is forced necessary by $\P_i$, and consider the
partial order $P_1\cross\cdots\cross\P_n$. The resulting forcing
extension $M^{\P_1\cross\cdots\cross\P_n}$ is a forcing extension
of each $M^{\P_i}$, and so each $\varphi_i$ holds there. Thus, any
given finite collection of sentences forceably necessary over a
fixed model of set theory is consistent. And so the whole
collection is consistent.\QED

\SubLemma. If $\varphi$ is forceably necessary, then this is
necessary; and if $\varphi$ is not forceably necessary, then this
also is necessary. In short, the collection of forceably
necessary sentences is invariant by
forcing.\label{NotForceablyNecessary}

\Proof: Suppose that $\varphi$ is forced necessary by the forcing
$\P$, so that $\varphi$ holds in every forcing extension of
$V^\P$. Given any other forcing extension $V^\Q$, one can still
force with $\P$ to obtain $V^{\Q\cross\P}$. Since this is the
same as $V^{\P\cross\Q}$, it is a forcing extension of $V^\P$.
Thus, over $V^\Q$ one could still force $\varphi$ to be
necessary; and so $\varphi$ is forceably necessary in $V^\Q$, and
hence necessarily forceably necessary in $V$. Conversely, suppose
that $\varphi$ is not forceably necessary. This means that
$\neg\varphi$ is forceable over every forcing extension. Thus,
over no forcing extension can $\varphi$ be forced necessary, so
$\varphi$ is necessarily not forceably necessary. In summary,
whether a statement is forceably necessary or not cannot be
affected by forcing.\QED

To prove the theorem, now, suppose that $M$ is any model of \ZFC.
Let $T$ be the collection of sentences that are forceably
necessary in $M$. This includes every axiom of \ZFC, since these
axioms hold in every forcing extension of $M$. By Lemma
\ref{Consistent}, the theory $T$ is consistent, and so it has a
model $N$. Suppose now that $\varphi$ is forceably necessary over
$N$; I aim to show that $\varphi$ is true in $N$. First, I claim
that $\varphi$ is actually in $T$. If not, then it is not
forceably necessary in $M$, and so by Lemma
\ref{NotForceablyNecessary} the sentence
$\neg\possible\necessary\varphi$ asserting that $\varphi$ is not
forceably necessary is necessary and hence in $T$. In this case,
since $N$ is a model of $T$, the sentence $\varphi$ cannot be
forceably necessary in $N$, contrary to our assumption. Thus,
$\varphi$ must be in $T$, and so it holds in $N$, as desired.\QED

\Proof A Second Proof: I would like now to give an alternative
proof of Theorem \ref{MP}, relying on a more traditional iterated
forcing construction. It is this iterated forcing argument that
will generalize to the case of $\MPtilde$.

Let me motivate the first lemma by mentioning that if there is an
inaccessible cardinal $\kappa$, then by the proof of the downward
Lowenheim-Skolem Theorem, one can find a closed unbounded set
$C\of\kappa$ of cardinals $\delta$ with $V_\delta\elesub
V_\kappa$. The structure $V_\kappa$, therefore, is a model of
\ZFC\ with a cardinal $\delta$ such that $V_\delta$ is an
elementary substructure of the universe. Axiomatizing this
situation, let the expression ``$V_\delta\elesub V$'' represent
the scheme, in the language with an additional constant symbol
for $\delta$, which asserts of any statement $\varphi$ in the
language of set theory that
$$\hbox{for every $x\in V_\delta$, if
$V_\delta\satisfies\varphi[x]$, then $\varphi(x)$.}$$ Each such
assertion in this scheme is first order, since one need only refer
to satisfaction in a set structure, $V_\delta$, and this is
provided by Tarski's definition of the satisfaction relation. The
little argument just given shows that if there is an inaccessible
cardinal, then there is a model of \ZFC\ satisfying the scheme
$V_\delta\elesub V$; indeed, since the club $C$ provides whole
towers of such $\delta$, the consistency strength of
$V_\delta\elesub V$ is easily seen to be strictly less than the
consistency of the existence of an inaccessible cardinal.

The really amazing thing, however, is that in fact
$\ZFC+V_\delta\elesub V$ is equiconsistent with \ZFC. One might
incorrectly guess that if \ZFC\ holds in $V$ and the scheme
$V_\delta\elesub V$ holds, then $V$ knows that $V_\delta$ is a
model of \ZFC; but this conclusion would confuse the `external'
\ZFC\ with the `internal' \ZFC\ of the model. What follows is only
that $V_\delta$ satisfies any particular instance of an axiom of
\ZFC\, but not the formula asserting that $V_\delta$ satisfies the
entire scheme \ZFC. This subtle distinction is crucial, in view of
the following elementary fact.\footnote{Thanks to Ali Enayat for
a helpful discussion of these ideas.}

\SubLemma. If\/ \ZFC\ is consistent, then so is
$\ZFC+V_\delta\elesub V$.\label{Elementary}

\Proof: Assume that \ZFC\ is consistent; so it has a model $M$. By
the L\'evy Reflection Theorem, every finite subcollection of the
theory $\ZFC+{V_\delta\elesub V}$ is modeled in some rank initial
segment of $M$, and therefore is consistent. So the whole theory
is consistent.\QED

% This is the old argument, which I replaced with the direct
% argument above:
%
% A direct argument via the Reflection Theorem is possible, but to
% avoid some confusing metamathematical issues, let me argue
% model-theoretically. If \ZFC\ is consistent, then standard
% model-theoretic arguments show that there is a countable model
% $M$ of \ZFC, which is {\df recursively saturated}; that is, every
% recursive type in finitely many parameters from $M$ that is
% finitely realized in $M$ is realized in $M$. Let $\Sigma$ be the
% scheme asserting $V_\delta\elesub V$, interpreted as a type in
% the variable $\delta$. Since the formulas in $\Sigma$, displayed
% above, are easily recognized, this is a recursive type. Further,
% it is finitely realized in any model of \ZFC\ because for any
% finite collection of formulas, the L\'evy Reflection Theorem
% provides an ordinal $\delta$ for which $V_\delta$ is absolute to
% $V$ for those formulas. Thus, by the recursive saturation of $M$,
% the type $\Sigma$ is realized in $M$. That is, there is an
% interpretation for $\delta$ in $M$ for which $V_\delta\elesub V$
% is true in $M$.\QED
%
% It is interesting and not difficult to notice that any $\delta$
% realizing the type $\Sigma$ in $M$ will have the feature that
% $(V_\delta)^M$ is itself recursively saturated. It follows, by a
% back-and-forth argument, that $M$ is isomorphic to
% $(V_\delta)^M$. Thus, any recursively saturated model of set
% theory is isomorphic to an initial segment of itself. These
% models are never well-founded; indeed, they cannot even have the
% correct $\omega$.

The theorem is now a consequence of the following lemma.

\SubLemma. Assume that $V_\delta\elesub V$. Then there is a
forcing extension of\/ $V$, by forcing of size at most $\delta$,
in which the Maximality Principle \MP\ holds.\label{ForcingMP}

\Proof: Since the Maximality Principle asserts, in a sense, that
all possible switches have been turned permanently on that are
possible to turn permanently on, the idea behind the proof of this
lemma will be an iteration that forces all statements that are
possible to force in a permanent way. The trick is to handle the
meta-mathematical issues that arise on account of the
non-definability of what is true or forceable. It is in doing this
that the hypothesis $V_\delta\elesub V$ will be used.

Let $\<\varphi_n\st n\in\omega>$ enumerate all sentences in the
language of set theory. I will define a certain forcing notion
$\P$, an $\omega$-iteration, all of whose initial segments $\P_n$
will be elements of $V_\delta$. So, suppose recursively that the
iteration $\P_n$ has been defined up to stage $n$, and consider
the model $V_\delta^{\P_n}$. If $\varphi_n$ is forceably necessary
over this model, then choose a poset $\Q_n$ which forces it to be
necessary, and let $\P_{n+1}=\P_n*\Qdot_n$ using a suitable name
for $\Q_n$. Let $\P$ be the finite support iteration of the
$\P_n$.

I claim that if $G\of\P$ is $V$-generic, then the Maximality
Principle holds in $V[G]$. To see this, suppose that $\varphi$ is
forceably necessary over $V[G]$. Necessarily, $\varphi$ is
$\varphi_n$ for some $n$. Factor the forcing $\P$ at stage $n$ as
$\P\iso\P_n*\Ptail$, where $\Ptail$ is the iteration after stage
$n$. Since $\varphi$ is forceably necessary over
$V[G]=V[G_n][\Gtail]$, it is also forceably necessary over
$V[G_n]$, because $V[G]$ is a forcing extension of $V[G_n]$ and
so any forcing extension of $V[G]$ is also a forcing extension of
$V[G_n]$. Since $V_\delta\elesub V$ and $\P_n\in V_\delta$,
whether a condition in $\P_n$ forces a given statement has the
same answer in $V$ as in $V_\delta$. Consequently,
$V_\delta[G_n]\elesub V[G_n]$. Thus, $\varphi_n$ is forceably
necessary over $V_\delta[G_n]$. In this case, the stage $n$
forcing $\Q_n$ forced it to be necessary, so $\varphi$ holds in
$V_\delta[G_{n+1}]$ and all its forcing extensions. By
elementarity again, $\varphi$ holds in $V[G_{n+1}]$ and all
forcing extensions. In particular, since $V[G]$ is a forcing
extension of $V[G_{n+1}]$, it holds in $V[G]$ and all further
extensions. That is, $\varphi$ is necessary in $V[G]$.\QED

This completes the second proof of Theorem \ref{MP}.\QED

\Corollary. The theory $\ZFC+\MP$ is equiconsistent with \ZFC.

Though the hypothesis $V_\delta\elesub V$ has now dropped away,
there is nevertheless a stronger connection between $\ZFC+\MP$ and
$\ZFC+V_\delta\elesub V$ than the previous corollary might
suggest, a connection I would like to explore in the next few
theorems. One cannot, for example, omit the hypothesis that
$V_\delta\elesub V$ from Lemma \ref{ForcingMP}.

\Theorem. If\/ \ZFC\ is consistent, then there is a model of\/
\ZFC\ having no forcing extension, and indeed, no extension of any
kind with the same ordinals, that is a model of\/ \MP. Thus, one
cannot omit the hypothesis that $V_\delta\elesub V$ from Lemma
\ref{ForcingMP}.

\Proof: I claim, first, that if there is a model of \ZFC, then
there is one in which the definable ordinals---and I mean ordinals
that are definable without parameters---are unbounded. To see
this, suppose that $W$ is a model of $\ZFC$, not necessarily
transitive, and let $M$ be the cut determined by the definable
elements of $W$. That is, $M$ consists of those elements in
$(V_\alpha)^W$, where $\alpha$ ranges over the definable ordinals
of $W$. I will show that $M\elesub W$, by verifying the
Tarski-Vaught criterion: suppose $W\satisfies \exists
x\,\varphi(x,a)$ for some $a\in M$. By the definition of $M$, it
must be that $a$ is in some $V_\alpha$ in the sense of $W$ for
some definable ordinal $\alpha$ of $W$. In $W$, let $\beta$ be the
least ordinal such that for every $y\in V_\alpha$ there is an
$x\in V_\beta$ such that $\varphi(x,y)$, if there is any such $x$
in $W$ at all. Such a $\beta$ exists by the Replacement Axiom.
Moreover, $\beta$ is definable in $W$, and consequently,
$(V_\beta)^W$ is included in $M$, and so $M$ has all the
witnesses it needs to verify the Tarski-Vaught criterion. Thus,
$M\elesub W$, and in particular, $M$ is a model of $\ZFC$ whose
definable elements are unbounded. By beginning with a model $W$
that also satisfies $V=L$, we may assume additionally that $M$
satisfies $V=L$ as well.

Suppose now, towards contradiction, that $M$ has an extension $N$,
with the same ordinals as $M$, that satisfies $\ZFC+\MP$.
Necessarily, $M=(L)^N$. Consider any definable ordinal $\alpha$ of
$W$, defined by the formula $\psi$. Since $M\elesub W$, the
formula $\psi$ also defines $\alpha$ in $M$, and since $M$ is the
$L$ of $N$, we know further that $\alpha$ is definable in $N$ as
``the ordinal satisfying $\psi$ in $L$''. Let $\varphi$ be the
sentence expressing, ``the ordinal defined by $\psi$ in $L$ is
countable''. This sentence is forceably necessary over $N$, since
if $\alpha$ is not countable in $N$, one can force to make it
countable, and having done this, of course, it remains countable
in any further extension. Thus, by $\MP$, the sentence $\varphi$
must be already true in $N$, and so $\alpha$ is countable in $N$.

The key observation is now that since the definable ordinals
$\alpha$ of $M$ were unbounded in $M$, it must be that every
ordinal of $M$ is countable in $N$, a contradiction, since these
two models of \ZFC\ have the same ordinals.\QED

While Lemma \ref{Elementary} establishes the equiconsistency of
\ZFC\ with $\ZFC+V_\delta\elesub V$, when it comes to transitive
models the two theories have different strengths. Specifically, if
$M$ is an $\omega$-model of $V_\delta\elesub V$, then $M$ has the
correct internal \ZFC, and so since every axiom of \ZFC\ is true
in $V_\delta$, the model $M$ agrees that $V_\delta$ is a model of
\ZFC. Consequently, inside any $\in$-minimal transitive model of
$\ZFC+V_\delta\elesub V$ there is a transitive model of \ZFC\ but,
by minimality, no transitive model of $\ZFC+V_\delta\elesub V$. So
it is relatively consistent that there is a transitive model of
the one theory, but not the other. In this sense, we might say
that the theory $\ZFC+V_\delta\elesub V$ has a greater
transitive-model consistency strength than \ZFC. The next theorem
shows that the transitive-model consistency strength of $\ZFC+\MP$
is similarly greater than $\ZFC$, since in fact it is the same as
$\ZFC+V_\delta\elesub V$.

\Theorem. The theories $\ZFC+\MP$ and $\ZFC+V_\delta\elesub V$ are
transitive-model-equiconsistent in the sense that one of them has
a transitive model if and only if the other also does. Indeed,
for any transitive model of one of these theories, there is a
model of the other with the same ordinals.

\Proof: Lemma \ref{ForcingMP} shows that every model of
$\ZFC+V_\delta\elesub V$ has a forcing extension that is a model
of $\ZFC+\MP$.

Conversely, suppose that $M$ is a transitive model of $\ZFC+\MP$.
Consider the $L$ of $M$, which must have the form $L_\gamma$, for
$\gamma=\ORD^M$. The arguments of the previous theorem establish
that every ordinal that is definable in $L^M$ is countable in $M$.
Let $\delta$ be the supremum of such ordinals. Thus, $\delta$ is
at most $\omega_1^M$, and in particular, $\delta<\gamma$. The
initial claim of the previous theorem establishes that in the
model $L^M$, the cut determined by the definable ordinals is an
elementary substructure. In this case, this means that
$V_\delta\elesub V$ in $L^M$, as desired.\QED

The argument establishes that for an arbitrary order-type $\xi$
(not necessarily well-ordered), the $\xi$-consistency strength of
$\ZFC+V_\delta\elesub V$ is at least as great as that of
$\ZFC+\MP$, since from every model of the former theory we can
build a model of the latter theory with the same ordinals. The
argument, however, does not quite establish the converse, because
in the non-transitive case one doesn't seem to know that the
supremum of the definable ordinals of $L^M$ exists as an ordinal
of $M$.

In the proof of Lemma \ref{ForcingMP}, the forcing $\P$ may have
size $\delta$, and in this case it may happen that every cardinal
below $\delta$ may be collapsed to $\omega$, so there is little
reason to expect that the axiom $V_\delta\elesub V$ remains true
in the forcing extension $V[G]$ (that is, using $(V[G])_\delta$
and $V[G]$). By taking a little care, however, the following
theorem explains how to arrange $\ZFC+\MP+V_\delta\elesub V$ in
the extension.

\Theorem. If\/ \ZFC\ is consistent, then so is the theory
${\ZFC}+{\MP}+{V_\delta\elesub V}$.

\Proof: First, one shows that if there is a model of \ZFC, then
there is a model of $\ZFC+{V_\delta\elesub V}+\cof(\delta)$ is
uncountable. This can be proved in the same way as Lemma
\ref{Elementary}, by simply adding the assertion
``$\cof(\delta)>\omega$'' to the theory. The L\'evy Reflection
Theorem shows that any finite collection of formulas is absolute
from $V$ to $V_\delta$ for a closed unbounded class of cardinals,
and so one may find such a $\delta$ of uncountable cofinality.

Given a model of $\ZFC+V_\delta\elesub V$ in which $\cof(\delta)$
is uncountable, one now proceeds with the argument of Lemma
\ref{ForcingMP}, using finite support in the iteration $\P$ of
that theorem. Since on cofinality grounds $\P$ has size less than
$\delta$, it follows that $V_\delta[G]\elesub V[G]$ in that
argument. That is, the axiom ``$V_\delta\elesub V$'' remains true
in $V[G]$, together with $\ZFC+\MP$.\QED

\Section Consistency of the Maximality Principle \MPtilde

Let me now turn to the boldface Maximality Principle \MPtilde, in
which arbitrary real parameters are allowed to appear in the
formulas $\varphi$ of the scheme
$\possible\necessary\varphi\implies\necessary\varphi$. As one
might expect, the boldface principle \MPtilde\ turns out to have a
stronger consistency strength.

Define the expression ``$\ORD$ is Mahlo'' to be the scheme
asserting that any closed unbounded class $C\of\ORD$ that is
definable from parameters contains a regular cardinal. This is
also sometimes referred to as the L\'evy Scheme (e.g., see
\cite{LevyScheme}). Whenever $\kappa$ is a Mahlo cardinal, then
$V_\kappa\satisfies\ZFC+\ORD$ is Mahlo; so the consistency
strength of ``$\ORD$ is Mahlo'' is strictly less than the
existence of a Mahlo cardinal. The converse of this implication,
however, is not generally true, because if $\kappa$ is a Mahlo
cardinal, then there is a club $C\of\kappa$ consisting of
$\delta$ with $V_\delta\elesub V$. These $V_\delta$ will
therefore also model ``$\ZFC+\ORD$ is Mahlo'', and so there will
be many non-Mahlo cardinals with that property.

\Theorem. The following theories are equiconsistent:
\begin{enumerate}
\item $\ZFC$ plus the Maximality Principle \MPtilde.
\item $\ZFC$ plus $V_\delta\elesub V$ for an inaccessible cardinal
$\delta$.
\item $\ZFC$ plus $\ORD$ is Mahlo.
\end{enumerate}\label{MPtilde}

\Proof: I will prove the theorem with a sequence of lemmas.

\SubLemma. Every model of theory 1 has an inner model of theory
2. Specifically, if $\MPtilde$ holds, then $\delta=\omega_1$ is
inaccessible in $L$, and $L_\delta\elesub L$.\label{inner}

\Proof: Assume $\MPtilde$ and let $\delta=\omega_1$. I claim,
first, that $\delta$ is inaccessible to reals. To see this,
suppose $x$ is a real, and let $\varphi$ be the statement ``the
$\omega_1$ of $L[x]$ is countable''. Since this is forceably
necessary, by $\MPtilde$ it is true. Thus,
$\omega_1^{L[x]}<\delta$, as I claimed. It follows from this that
$\delta$ is inaccessible in $L$.

Second, I claim that $L_\delta$ is an elementary substructure of
$L$ (note: I prove this part of the Lemma only as a scheme). To
see this, simply verify the Tarski-Vaught criterion: suppose that
$L$ satisfies $\exists y\, \psi(a,y)$ for a set $a$ in $L_\delta$.
Let $\alpha$ be least such that there is such a $y$ in
$L_\alpha$. Consider the statement $\varphi$ asserting ``the
least $\alpha$, such that there is a $y$ in $L_\alpha$ with
$\psi(a,y)^L$, is countable". This is expressed using the
parameter $a$, which is coded with a real. Since it is forceably
necessary, it must be true by $\MPtilde$, and so $\alpha$ is
countable. Thus, $y$ can be found in $L_\delta$, and so the
Tarski-Vaught criterion holds. In summary, $L$ satisfies theory
2.\QED

\SubLemma. If theory 2 holds, then theory 1 holds in a forcing
extension. Indeed, if\/ $V_\delta\elesub V$ and $\delta$ is an
inaccessible cardinal, then after the L\'evy collapse making
$\delta$ the $\omega_1$ of the extension, the Maximality
Principle \MPtilde\ holds.\label{outer}

\Proof: Assume that $V_\delta\elesub V$ and $\delta$ is
inaccessible. In order to obtain \MPtilde\ in a forcing extension,
I will define a certain finite support $\delta$-iteration $\P$
which at every stage forces with some poset of rank less than
$\delta$. Let $\<\varphi_\alpha\st\alpha<\delta>$ enumerate, with
unbounded repetition, all possible sentences in the forcing
language for initial segments of such iterations with parameters
having names in $V_\delta$. Suppose now that the iteration
$\P_\alpha$ is defined up to stage $\alpha$; I would like to
define the stage $\alpha$ forcing $\dot\Q_\alpha$. If
$\varphi_\alpha$ is a sentence with parameters naming objects in
$V_\delta^{\P_\alpha}$, and if it is possible to force over
$V_\delta^{\P_\alpha}$ so as to make $\varphi_\alpha$ necessary,
then choose $\dot\Q_\alpha\in V_\delta$ to be (the name of) a
poset accomplishing this. Otherwise, let $\dot\Q_\alpha$ be
trivial forcing. This defines the $\delta$-iteration $\P$.

Suppose that $G\of\P$ is $V$-generic, and consider the model
$V[G]$. First observe that the cardinal $\delta$ is not collapsed
by this forcing, since it is a finite support iteration of
$\delta$-c.c. forcing. Consequently, every real in $V[G]$ appears
in some $V[G_\alpha]$. Thus, if $\varphi$ is a sentence with
(names for) real parameters from $V[G]$, it must be
$\varphi_\alpha$ for some $\alpha$, with the names in $\varphi$
coming from $V[G_\alpha]$. So it suffices to consider the formulas
$\varphi_\alpha$.

Suppose that such a $\varphi_\alpha$ is forceably necessary over
$V[G]$. Then it is also forceably necessary over $V[G_\alpha]$,
since $V[G]$ is a forcing extension of $V[G_\alpha]$. Hence, by
the elementarity $V_\delta[G_\alpha]\elesub V[G_\alpha]$ (as in
Lemma \ref{ForcingMP}), the sentence $\varphi_\alpha$ is forceably
necessary over $V_\delta[G_\alpha]$. The stage $\alpha$ forcing
$\Q_\alpha$ must then have forced it to be necessary, and so
$\varphi_\alpha$ became necessary in $V_\delta[G_{\alpha+1}]$, and
hence also in $V[G_{\alpha+1}]$. Thus, since $V[G]$ is obtained by
further forcing over this model, $\varphi_\alpha$ remains
necessary in $V[G]$. So $V[G]$ satisfies \MPtilde.

Finally, let me consider the precise nature of the iteration
$\P$; in fact, we know $\P$ very well. Since the assertion ``$a$
is countable'' is always forceably necessary, it follows that
every element $a\in V_\delta$ is countable in $V[G]$. That is,
unboundedly often, the forcing $\Qdot_\alpha$ collapses elements
of $V_\delta$. Since $\delta$ itself is not collapsed, it becomes
the $\omega_1$ of $V[G]$. But up to isomorphism, the only finite
support $\delta$-iteration of forcing notions of size less than
$\delta$ that collapses all cardinals below $\delta$ to $\omega$
is the L\'evy collapse of $\delta$ to $\omega_1$. Thus, $\P$ is
isomorphic to the L\'evy collapse, as I claimed.\QED

David Asper\'o has used a similar iteration to prove a similar
conclusion relating to his generalized bounded forcing axioms
(see Theorem 3.7 in \cite{DavidAspero}).

\SubLemma. Theory 2 implies theory 3.

\Proof: Suppose that $V_\delta\elesub V$ and $\delta$ is
inaccessible. If $C\of\ORD$ is a definable club, then of course
$C\intersect \delta$ is unbounded in $\delta$. Thus, $\delta\in
C$, and so $C$ contains a regular cardinal, as desired.\QED

\SubLemma. If theory 3 is consistent, then so is theory 2.

\Proof: This argument is similar to Lemma \ref{Elementary}.
Suppose that $M$ is a model of \ZFC\ plus the scheme asserting
that $\ORD$ is Mahlo. By the L\'evy Reflection Theorem, any
finite collection of formulas reflects from $V$ to $V_\theta$ for
a closed unbounded class of cardinals $\theta$. Thus, since
$\ORD$ is Mahlo in $M$, we can find in $M$ a regular cardinal
$\delta$ such that the formulas are absolute between $V$ and
$V_\delta$ in $M$. We may assume that the finite collection of
formulas includes the assertion that for every $\beta$, the
cardinal $2^\beta$ exists, and from this it follows that $\delta$
is a strong limit cardinal. Thus, $\delta$ is actually
inaccessible in $M$. So every finite subcollection of theory 2 is
consistent, and so the entire theory is consistent.\QED

This completes the proof of the theorem.\QED

Lemmas \ref{inner} and \ref{outer} establish that if there is a
model of one of the two theories ``$\ZFC+\MPtilde$" or
``$\ZFC+{V_\delta\elesub V}+\delta$ is inaccessible," then there
is a model of the other with the same ordinals. Therefore, we may
conclude:

\Corollary. The theories ``$\ZFC+\MPtilde$'' and
``$\ZFC+{V_\delta\elesub V}+{\delta}$ is inaccessible'' are
$\xi$-equiconsistent for every order type $\xi$.

Since Lemma \ref{outer} proceeds via the L\'evy collapse, one
naturally obtains the \CH\ in the extension, and this leads one to
wonder whether $\MPtilde+\neg\CH$ is consistent. The corresponding
question does not arise with the lightface Maximality Principle
\MP\ because once \MP\ is true, it persists to all forcing
extensions, and so one can simply force \CH\ or $\neg\CH$ as
desired, while preserving \MP. But with the boldface Maximality
Principle, it is not clear that $\MPtilde$ persists to an
extension where new real parameters are available. The following
theorem shows how to sidestep this worry.

\Theorem. The boldface Maximality Principle \MPtilde\ is
relatively consistent with either \CH\ or $\neg\CH$.

\Proof: Over any model one can force \CH\ without adding reals,
and this will preserve \MPtilde. Alternatively, the model of
Lemma \ref{outer} satisfies $\MPtilde+\CH$, since it is obtained
via the L\'evy collapse of an inaccessible cardinal.

To obtain a model of $\MPtilde+\neg\CH$, we will simply add Cohen
reals over this model.\footnote{I am grateful to Ilijas Farah for
pointing out the applicability of Solovay's technique here.}
Specifically, I claim that the model $V[G][g]$, where $G\of\P$ is
$V$-generic for the iteration of Lemma \ref{outer} and $g$ is
$V[G]$-generic for the forcing to add at least $\omega_2^{V[G]}$
many Cohen reals, is a model of $\MPtilde+\neg\CH$. Certainly
$\neg\CH$ is no problem, so consider \MPtilde. The instances of
the axiom scheme
$\possible\necessary\varphi\implies\necessary\varphi$, where
$\varphi$ has real parameters from $V[G]$, hold in $V[G]$ by
Lemma \ref{outer}, and each of these instances persists to every
forcing extension, including $V[G][g]$. At issue is whether the
scheme holds of formulas using the new parameters available in
$V[G][g]$. Consider, therefore, a real $x$ in $V[G][g]$ and the
corresponding extension $V[G][x]$. Since $x$ is added by a
countable forcing notion $\Q$, we may reorganize this forcing as
$V[G][x]=V[x][G']$, where $x$ is now $V$-generic for some
countable forcing notion and $G'$ is $V[x]$-generic for the
L\'evy collapse of $\kappa$ in $V[x]$ (this is because the
quotient forcing $(\P*\Q)/x$ has size $\kappa$, is $\kappa$-c.c.
and collapses every cardinal below $\kappa$ to $\omega$; hence it
is the L\'evy collapse). Since $V_\kappa[x]\elesub V[x]$, it
follows by Lemma \ref{outer} that \MPtilde\ holds in $V[x][G']$.
Thus, instances of the \MPtilde\ scheme using the parameter $x$
hold in $V[x][G']=V[G][x]$. Since these persist to any forcing
extension, they also hold in $V[G][g]$. So \MPtilde\ holds in
$V[G][g]$, as desired.\QED

Techniques of Kunen, employing the same reorganization idea as
above, can be used to obtain a model of $\MPtilde+\neg\CH+\MA$.

Let me close this section with a discussion of the relative
consistency of \MP\ with large cardinals.

\Theorem. The Maximality Principles \MP\ and \MPtilde\ are
consistent with the existence of inaccessible cardinals,
measurable cardinals, supercompact cardinals, or what have you,
assuming the consistency of a proper class of such cardinals
initially (plus, for \MPtilde, that $\ORD$ is Mahlo). On the other
hand, the Maximality Principle is also consistent with the
nonexistence of such cardinals, assuming only the consistency of
\ZFC.

\Proof: Suppose that $\ZFC$ is consistent with the existence of a
proper class of, for example, measurable cardinals. The argument
of Lemma \ref{Elementary} produces a model of
$\ZFC+V_\delta\elesub V$ with a proper class of measurable
cardinals. Since the forcing of Lemma \ref{ForcingMP} has size at
most $\delta$, the measurable cardinals above $\delta$ survive to
the forcing extension in which $\MP$ holds. For $\MPtilde$, one
should add the hypothesis that ``$\ORD$ is Mahlo'' to this
argument in order to get a model as above for which $\delta$ is
inaccessible, and then carry on as in Theorem \ref{MPtilde}. This
procedure shows that any large cardinal that is not destroyed by
small forcing will be consistent with \MP\ or \MPtilde.

To obtain the Maximality Principle in the absence of all large
cardinals, one should simply apply the argument of Theorem
\ref{MP} using a model with no inaccessible cardinals. That is,
if $M$ is any model of \ZFC\ with no inaccessible cardinals, then
no forcing extension of $M$ will have inaccessible cardinals, and
so the assertion that there are no inaccessible cardinals will be
necessary in $M$ and hence in the theory $T$ of that proof.\QED

The attentive reader will observe that in order to deduce the
consistency of the Maximality Principle with an inaccessible
cardinal, the previous proof begins with and eventually produces
models with proper classes of such cardinals. This apparent
inefficiency suggests the question of whether the Maximality
Principles are consistent with the existence of a single
inaccessible cardinal. The answer, perhaps surprising at first, is
that they are not.

\Theorem. If the Maximality Principle \MP\ holds, then there are
a proper class of inaccessible cardinals, if
any.\footnote{Chalons \cite{Chalons2} uses a similar idea when
proving that the Chalons Maximality Principle does not follow
from any large cardinal axiom.}\label{IfAny}

\Proof: Suppose that the Maximality Principle holds, but that the
inaccessible cardinals are bounded. Then the assertion ``there are
no inaccessible cardinals'' is forceably necessary, because one
could force, in a permanent way, to make all the inaccessible
cardinals countable. Thus, under \MP, there must have been no
inaccessible cardinals to begin with.\QED

The same argument works with any other large cardinal whose
existence is downward absolute from a forcing extension, such as
Mahlo cardinals. For measurable cardinals, whose existence is not
downward absolute in this way, the situation is slightly more
complicated. But one can define that a cardinal $\kappa$ is {\df
potentially} measurable when it is measurable in a forcing
extension. This notion, of course, is downward absolute from any
forcing extension. Consequently, if \MP\ holds and there is a
potentially measurable cardinal, then there must be a proper class
of potentially measurable cardinals. (If there are a bounded
number, then force to collapse them all, making none, in a
permanent way; so there must have been none to begin with.) An
identical fact holds for potentially supercompact cardinals, or
what have you.

\Section The Necessary Maximality Principle

We turn now to the strongest principle I have mentioned in this
paper, the Necessary Maximality Principle $\necessary\MPtilde$. I
will begin with a theorem hinting at surprising strength.

\Theorem. If\/ $\necessary\MPtilde$ holds, then $0^\#$ exists.

\Proof: It suffices, by the Covering Lemma of Jensen, to find a
violation of covering. For this, it suffices to show that $L$
fails to compute the successors of singular cardinals correctly.
In fact, I will show that $L$ fails to compute the successor of
any cardinal correctly. Consider any cardinal $\kappa$ and its
successors $(\kappa^\plus)^L$ and $(\kappa^\plus)^V$ in the two
models. Let $g$ be $V$-generic for the canonical forcing that
collapses $\kappa$ to $\omega$. This forcing is the same in $V$ or
$L$ and has size $\kappa$ in either model, and so by the chain
condition it preserves $\kappa^\plus$ over either model, making it
the $\omega_1$ of the extension. The key observation now is that
the sentence $\varphi$, asserting ``$(\omega_1)^{L[g]}$ is
countable,'' is forceably necessary. Thus, by \MPtilde\ in $V[g]$,
this sentence is true in $V[g]$. Putting these conclusions
together, we have
$$(\kappa^\plus)^L=(\omega_1)^{L[g]}<(\omega_1)^{V[g]}=(\kappa^\plus)^V,$$
and so $L$ does not compute the successor of $\kappa$
correctly.\QED

This argument generalizes to the following, where by
$W_\gamma$ I mean $(V_\gamma)^W$.

\Theorem. Suppose that\/ $\necessary\MPtilde$ holds and $W$ is a
definable transitive class, invariant under set forcing. Then:
\begin{enumerate}
\item The class $W$ does not compute successor cardinals correctly.
\item Indeed, every regular cardinal is inaccessible in $W$.
\item Indeed, every cardinal is a limit of inaccessible cardinals
in $W$.
\item For every cardinal $\gamma$, $W_\gamma\elesub W$.
\item For every cardinal $\gamma$, $H(\gamma)^W\elesub W$.
\item For every cardinal $\gamma$, if a set $x$ is definable in
$W$ from an object $a\in H(\gamma)^V$, then $x\in H(\gamma)^V$.
\end{enumerate}\label{W}

\Proof: Clearly, statement $1$ of the theorem is implied by
statement $2$, which is in turn implied by statement $3$. Let me
argue that statement $3$ is implied by statement $4$. If
$W_\gamma\elesub W$, then clearly $\gamma$ is a strong limit in
$W$, and hence, if regular, $\gamma$ is inaccessible in $W$.
Thus, the inaccessible cardinals are unbounded in $W$, and hence
also in $W_\gamma$, so statement $3$ holds. Statement $4$ easily
follows from statement $5$, because if $H(\gamma)^W\elesub W$,
then $H(\gamma)^W$ has the correct $(V_\alpha)^W$ for
$\alpha<\gamma$, and so $H(\gamma)^W=W_\gamma$, giving statement
$4$.

Next, let me show that statement $5$ follows from statement $6$,
by verifying the Tarski-Vaught criterion. Suppose that $W$
satisfies $\exists x\,\psi(x,a)$ where $a\in H(\gamma)^W$. The
least hereditary size in $W$ of such an $x$ is definable in $W$
from $a$, and so by statement $6$ this size must be less than
$\gamma$. In particular, there will be such an $x$ in
$H(\gamma)^W$, and so the criterion is fulfilled.

It remains only to prove statement $6$ of the theorem. Suppose
that $x$ is definable in $W$ from the object $a\in H(\gamma)$.
Let $\beta<\gamma$ be the hereditary size of $a$, so that $a\in
H(\beta^\plus)$, and suppose that $g\of\beta$ is $V$-generic for
the canonical forcing that collapses $\beta$ to $\omega$. Let
$\varphi$ be the sentence asserting ``$x$ is hereditarily
countable'', using the definition of $x$ in $W$ from $a$, coded
with the real $g$. Thus, $\varphi$ is a sentence with real
parameters, and it is forceably necessary (please note that I am
using here the invariance of $W$ under under set forcing in order
to know that $x$ is still defined in the same way in any forcing
extension). Consequently, since $\necessary\MPtilde$ implies
$\MPtilde$ in $V[g]$, the sentence $\varphi$ must be true in
$V[g]$. That is, $x$ is hereditarily countable in $V[g]$. Since
$\beta^\plus$ and $\gamma$ are preserved to $V[g]$, it follows
that $x$ has hereditary size less than $\gamma$ in $V$, as
desired. So the proof is complete.\QED

The theorem is easily modified to handle the situation when $W$ is
defined from parameters. Indeed, if $W$ is definable from a real
parameter, the previous proof goes through essentially
unchanged. More generally, we have:

\Corollary. Suppose that\/ $\necessary\MPtilde$ holds and $W$ is a
transitive class, definable from parameters $\vec p\in H(\eta^\plus)$,
which is invariant under set forcing. Then the conclusions of the
previous theorem hold above $\eta$.

\Proof: By this, I mean that such a class $W$ will not compute
successor cardinals above $\eta$ correctly, that every cardinal
above $\eta$ will be a limit of inaccessible cardinals in $W$, and
so on, and that every set $x$ definable in $W$ from an object in
$H(\gamma)$ for some $\gamma$ above $\eta$ will be in $H(\gamma)$.

Suppose, therefore, that $W$ is definable from the parameters
$\vec p\in H(\eta^\plus)$, and that it is invariant by set
forcing. In the forcing extension $V[g]$ that collapses $\eta$ to
$\omega$, the parameters $\vec p$ become hereditarily countable,
and the class $W$ therefore becomes definable from a real $z$.
Thus, applying the previous proof to the class $W$, defined by the
real $z$ in $V[g]$, the conclusion follows for the cardinals of
$V[g]$. Since these are precisely the cardinals above $\eta$ in
$V$, the corollary is proved.\QED

\Corollary. If\/ $\necessary\MPtilde$ holds, then for every set
$x$, the set $x^\#$ exists. That is, the universe is closed under
sharps.

\Proof: Fix any set $x$, and let $g$ be a real collapsing
cardinals so that $x$ can be coded by a real $\bar x$ in $V[g]$.
Let $W=L[x]$. This is definable from $\bar x$, and is invariant
in any further extension of $V[g]$. Thus, by the previous
corollary, if $\necessary\MPtilde$ holds, then $L[x]$ does not
compute the successors of singular cardinals above $x$ correctly.
It follows that $x^\#$ exists in $V[g]$, and so it also exists in
$V$.\QED

The same idea shows that the universe is closed under daggers and
pistols. Any model of the form $L[x,\mu]$, for example, where the
measure $\mu$ lives on the smallest possible ordinal, is definable
from $x$, and Theorem \ref{W} shows that this model cannot compute
successor cardinals correctly above $x$.

\Theorem Projective Absoluteness Theorem. The Necessary
Maximality Principle $\necessary\MPtilde$ implies Projective
Absoluteness. That is, if\/ $\necessary\MPtilde$ holds, then
projective truth is invariant by set forcing; every $\Sigma^1_n$
formula is absolute to any set forcing extension.

\Proof: Assume $\necessary\MPtilde$. I will show, by induction on
$n$, that every $\Sigma^1_n$ formula $\varphi$ is absolute to any
set forcing extension. This is trivial for $n=0$ (and indeed, it
holds automatically up to $\Sigma^1_2$ by the Shoenfield
Absoluteness Theorem). Assume by induction that it is true for all
$\Sigma^1_n$ formulas. Further, since $\necessary\MPtilde$ holds
in all set forcing extensions, suppose that this induction
hypothesis also holds in all set forcing extensions, that is,
that the $\Sigma^1_n$ formulas are absolute from any set forcing
extension to any further set forcing extension. Since the
collection of $\varphi$ for which this hypothesis is true is
clearly preserved under Boolean combinations, it suffices to
consider an existential formula $\psi(x)=\exists
y\,\varphi(x,y)$, where $\varphi$ is absolute from any forcing
extension to any further forcing extension. If $\psi(a)$ is true
in such a model, then since $\varphi$ is absolute to any set
forcing extension and both $a$ and the witness $y$ still exist
there, $\psi(a)$ remains true in any set forcing extension.
Conversely, suppose that $\psi(a)$ is true in a set forcing
extension $V^\P$, with $a\in V$. Thus, the witness $y$ to
$\varphi$ was added by the forcing $\P$. Since by the induction
hypothesis the formula $\varphi(a,y)$ is absolute to any further
forcing extension, the formula $\psi(a)$ is necessary in $V^\P$
and therefore forceably necessary in $V$. Thus, by
$\necessary\MPtilde$, it is already true in $V$, as desired. So
every projective formula is absolute by set forcing, and the
theorem is proved.\QED

\Corollary. The consistency of\/ $\necessary\MPtilde$ implies that
of the existence of infinitely many strong cardinals.

\Proof: Kai Hauser \cite{ProjAbs} has shown that projective
absoluteness is equiconsistent with the existence of infinitely
many strong cardinals.\QED

The next theorem hints at an even stronger strength for
$\necessary\MPtilde$, assuming that $\ORD$ is ineffible, i.e., the
scheme asserting that for every class partition $F:[\ORD]^2\to 2$
there is a homogeneous stationary class $H\of\ORD$. This axiom is
naturally expressed in G\"odel-Bernays set theory, using class
quantifiers. I am indebted to Philip Welch for his helpful remarks
concerning the next few corollaries.

\Corollary. If\/ $\necessary\MPtilde$ holds and $\ORD$ is
ineffable, then there is an inner model of a Woodin cardinal.

\Proof: Steel \cite{Steel} has proved that if there is no inner
model of a Woodin cardinal and if the universe is closed under
sharps and $\ORD$ is ineffable, then there is a definable class
$K$, now widely known as the Mitchell-Steel core model, which is
invariant under set forcing and which computes the successors of
singular cardinals correctly. Since Theorem \ref{W} shows that if
$\necessary\MPtilde$ holds, then the universe is closed under
sharps but there can be no such class $K$, it follows under the
hypothesis that there must be an inner model with a Woodin
cardinal.\QED

One can modify the assumption that $\ORD$ is ineffable by building
the core model up to a weakly compact cardinal $\kappa$.
Specifically, Steel \cite{Steel} shows that if $\kappa$ is weakly
compact, $V_\kappa$ is closed under sharps and $V_\kappa$ has no
inner model of a Woodin cardinal, then the core model $K$ built
up to $\kappa$ will compute the successors of singular cardinals
correctly. Thus, since $\MPtilde$ implies closure under sharps
and refutes the possibility of such a class $K$, this establishes:

\Corollary. If $\kappa$ is weakly compact and
$V_\kappa\satisfies\necessary\MPtilde$, then $V_\kappa$ has an
inner model of a Woodin cardinal.

This kind of analysis can be lifted to those core models with more
Woodin cardinals, provided that Weak Covering holds. Following
this line, Philip Welch has announced the following:

\Theorem. (Welch) If\/ $\necessary\MPtilde$, then Projective
Determinacy \PD\ holds.

Welch also conjectured that if $\necessary\MPtilde$ holds and
there is an inaccessible cardinal, then $\AD^{L(\R)}$. Woodin has
announced, in conversation, that $\necessary\MPtilde$ implies
$\AD^{L(\R)}$ outright.

I regret that the obvious problem remains unsolved:

\Theorem Challenge Problem. Is $\ZFC+\necessary\MPtilde$
(relatively) consistent?

I would welcome a consistency proof relative to any large
cardinal hypothesis.

Let me present a theorem, suggested by W. Hugh Woodin, that makes
at least some progress in this direction. Let's say that a
sentence is {\df local} if it can be expressed in the form ``{\it
$\psi$ holds in $V_\kappa$, where $\kappa$ is the least
inaccessible cardinal}''. Such sentences assert a truth that can
be checked locally, in this $V_\kappa$, without need to consult
the entirety of $V$. This is a broad class of assertions,
including all the conjectures and theorems of classical
mathematics, as well as the bulk of contemporary mathematics.

Please note that the least inaccessible cardinal, of course, is
not invariant by forcing; a model may disagree with a forcing
extension on the particular value of $\kappa$ even when the two
models agree about some particular local assertion.

\Theorem. If there are a proper class of Woodin cardinals, then
the Necessary Maximality Principle $\necessary\MPtilde$ holds for
local assertions (using real parameters in any forcing extension).

\Proof: I will first prove that if there are a proper class of
Woodin cardinals, then the Maximality Principle $\MPtilde$ holds
for local assertions using real parameters in the ground model.
From this, since the hypothesis of a proper class of Woodin
cardinals is preserved by set forcing, it follows that the same
conclusion holds in every forcing extension as well. That is, from
the hypothesis of a proper class of Woodin cardinals one can
conclude the Necessary Maximality Principle $\necessary\MPtilde$
for local assertions as well.

To begin, then, suppose that $\varphi$ is local and forceably
necessary, forced necessary by the forcing $\P$, with real
parameters in $V$. Choose a Woodin cardinal $\delta$ above the
rank of $\P$ and let $\P_{<\delta}$ be the stationary tower
forcing corresponding to $\delta$. Suppose that
$G_{<\delta}\of\P_{<\delta}$ is $V$-generic below a condition
forcing that $2^\card{\P}$ becomes countable. By the basic
properties of stationary tower forcing (see, e.g.
\cite{StationaryTower}, \cite{StationaryTower2}), there is a
generic elementary embedding $j:V\to M\of V[G_{<\delta}]$, with
$M^{<\delta}\of M$ in $V[G_{<\delta}]$. The model $M$ therefore
agrees with $V[G_{<\delta}]$ well beyond the first inaccessible
cardinal. Since furthermore $2^\card{\P}$ was collapsed, in
$V[G_{<\delta}]$ we may construct a $V$-generic filter for $\P$,
and therefore we may view $V[G_{<\delta}]$ as a forcing extension
of $V^\P$. Since $\varphi$ is made necessary in $V^\P$, it
follows that $\varphi$ holds in $V[G_{<\delta}]$. Thus, since
$V[G_{<\delta}]$ agrees with $M$ well beyond the first
inaccessible cardinal, $\varphi$ holds also in $M$. Finally,
since the real parameters of $\varphi$ are fixed by $j$, it
follows by elementarity that $\varphi$ holds also in $V$, as
desired.\QED

In the previous argument, one can clearly relax the notion of
``local'' assertions quite a bit. Rather than using the least
inaccessible cardinal, for example, one could use the least Mahlo
cardinal or the second measurable cardinal or the fifth measurable
limit of measurable cardinals, and so on. All that is required in
the argument is that the cardinal in question be less than or
equal to the Woodin cardinal $\delta$ used in the proof, and that
the cardinal be the same in $M$ and $V[G_\delta]$.

Finally, there is a seeming affinity between the Maximality
Principles and the possibility that the theory of $L(R)$ is
persistent, that is, that it is invariant by set forcing, as it
is if there are sufficiently many Woodin cardinals. But what is
the exact connection between these notions?

\Question. What is the relationship between $\necessary\MPtilde$
and the $\Th(L(\R))$ being invariant by forcing?

\Section Modified Maximality Principles

The Maximality Principle admits numerous modified forms, obtained
by restricting the kinds of forcing extensions allowed in the
modal quantifiers. Specifically, suppose that $\cal P$ is a
definable class of forcing notions; we may define that a
statement $\varphi$ is possible or forceable {\df by forcing
notions in $\cal P$}, written $\possible_{\cal P} \varphi$, when
it holds in a forcing extension $V^\P$ by a forcing notion in
$\cal P$. Similarly, $\varphi$ is necessary or persistent {\df by
forcing notions in $\cal P$}, written $\necessary_{\cal P}
\varphi$, when $\varphi$ holds in $V$ and all extension $V^\P$ by
forcing $\P\in\cal P$. Combining the two notions, we see that
$\varphi$ is forceably necessary by forcing notions in $\cal P$,
written $\possible_{\cal P}\necessary_{\cal P}\varphi$, when it
holds in some extension $V^\P$ and all subsequent extensions
$V^{\P*\Qdot}$, where $\P$ and $\Qdot$ are taken from $\cal P$
(interpreting $\cal P$, in the case of $\Qdot$, in the model
$V^\P$).

The Maximality Principle can now be stated in the restricted form
as follows:

\QuietTheorem Maximality Principle ($\MP_{\cal P}$). If a
statement is forceably necessary by forcing notions in $\cal P$,
then it is persistent by forcing in $\cal P$.

\noindent Natural instances of this principle include:

\begin{itemize}
\item ($\MP_{card}$)
Any statement forceably necessary by forcing that preserves
cardinals is persistent by such forcing.

\item ($\MP_{\omega_1}$) Any
statement forceably necessary by forcing that preserves $\omega_1$
is persistent by such forcing.

\item ($\MP_{ccc}$) Any statement
forceably necessary by c.c.c.~forcing is persistent by
c.c.c.~forcing.

\item ($\MP_{proper}$) Any statement forceably necessary by
proper forcing is persistent by such forcing.

\item ($\MP_{semi-proper}$) Any statement forceably necessary by
semi-proper forcing is persistent by such forcing.

\item ($\MP_{\hbox{\scriptsize\sc ch}}$) Any statement forceably necessary by forcing
preserving \CH\ is persistent by such forcing.

\end{itemize}

And so on. All these principles admit the equivalent forms of
Theorem \ref{EquivalentForms}, so that $\MP_{\hbox{\scriptsize\sc
ch}}$, for example, is equivalent to the assertion that every
true sentence is forceable by \CH-preserving forcing over every
forcing extension with \CH.

But are they consistent? One might initially think, since these
principles simply restrict the class of forcing notions in the
Maximality Principle, that they are special cases of or at least
follow immediately from the full Maximality Principle \MP. But
this is not correct. For the reader will readily realize that an
assertion may be forceably necessary by c.c.c.~forcing, for
example, but not be forceably necessary at all (e.g. $\neg\CH$).
An assertion can be persistent under c.c.c.~forcing, but not
under all forcing.

So it is not immediately clear that these modified Maximality
Principles are consistent. I will prove now that at least some of
them are.

\Theorem. Suppose that $\cal P$ is necessarily closed under
finite support or countable support iterations of countable
length. Assuming $V_\delta\elesub V$, then, there is a forcing
extension obtained by forcing with a notion in $\cal P$ in which
$\MP_{\cal P}$ holds.

\Proof: In the proof of Lemma \ref{ForcingMP}, let us simply
ensure that the forcing $\Q_n$ at stage $n$ comes from $\cal P$.
Then, define $\P$ to be the finite support or countable support
iteration of the $\Q_n$, depending on the closure of $\cal P$.
Having done so, the iteration $\P$, as well as the tail forcing
$\Ptail$, will be in $\cal P$, and so the argument of Lemma
\ref{ForcingMP} works to establish $\MP_{\cal P}$ in $V^\P$.\QED

\Corollary. The principles $\MP_{ccc}$, $\MP_{proper}$ and
$\MP_{semi-proper}$ are equiconsistent with \ZFC.

\Theorem. Assume that \ZFC\ is consistent. If in any model of set
theory, the class $\cal P$ is closed under finite iterations and
members of $\cal P$ are necessarily in $\cal P$, then
$\ZFC+\MP_{\cal P}$ is consistent.

\Proof: This argument follows the first proof of Theorem \ref{MP}.
Suppose that $M$ is any model of \ZFC, and let $T$ be the
collection of sentences that are forceably necessary over $M$,
referring throughout this argument only to forcing notions in
$\cal P$.

I claim first that the theory $T$ is consistent. Given any finite
subcollection $\varphi_1,\ldots,\varphi_n$, with $\varphi_i$
forced necessary by $\P_i$, consider
$\P=\P_1\cross\cdots\cross\P_n$. By the hypotheses of the
theorem, this is in $\cal P$. Further, the extension $M^\P$ is a
forcing extension of each $M^{\P_i}$, by forcing in $\cal P$
there. Thus, the extension $M^\P$ satisfies each $\varphi_i$, and
so this collection is consistent.

Next, I claim that $T$ is absolute to all forcing extensions of
$M$. If a sentence is forceably necessary over $M$ with respect to
forcing in $\cal P$, then the forcing that made it necessary will
still be in $\cal P$ in any forcing extension, and serve to make
it necessary there. So it will be forceably necessary over any
forcing extension also. Conversely, if a sentence is forceably
necessary over a forcing extension of $M$ by forcing in $\cal P$,
then clearly it is forceably necessary in $M$.

Finally, I claim that any model $N$ of $T$ is a model of
$\ZFC+\MP_{\cal P}$. Certainly $\ZFC\of T$ since these are
necessary over $M$. If $\varphi$ is forceably necessary over $N$,
then $\varphi$ must be in $T$, for if not, then the assertion
that $\varphi$ is not forceably necessary would be necessary and
hence in $T$, a contradiction. Thus, $\varphi\in T$ and so
$\varphi$ holds in $N$, as desired.\QED

Define that a poset is {\df necessarily-c.c.c.} if it is
c.c.c.~in every c.c.c.~forcing extension. For example, while
Souslin trees lose their c.c.c.~property after forcing with them,
the forcing to add (any number of) Cohen reals does not. Let
$\MP_{\necessary ccc}$ be the corresponding Maximality Principle
for such forcing notions. Similarly, define that $\P$ is {\df
necessarily $\omega_1$-preserving} if it preserves $\omega_1$
over every forcing extension with the correct $\omega_1$, and
denote the corresponding Maximality Principle by
$\MP_{\necessary\omega_1}$.

\Corollary. The principles $\MP_{\necessary ccc}$ and
$\MP_{\necessary\omega_1}$ are equiconsistent with \ZFC.

\Proof: The corresponding classes of posets have the closure
properties of the theorem.\QED

Let me turn now to boldface versions of the modified principles.
The intriguing possibility here is the potential to go beyond the
reals by including some uncountable parameters in the scheme. In
Observation \ref{Parameters}, the argument that excluded
uncountable parameters, we used forcing that collapsed cardinals.
So when restricting to a class of posets that preserve all
cardinals, such as the c.c.c.~posets, we might hope to include
many more parameters.

But unfortunately, we cannot hope to use any set as a parameter.
In the c.c.c.~case, for example, suppose that $A$ is the set of
reals of a model $M$. The assertion, ``{\it there is a real
$x\notin A$},'' is forceably necessary by c.c.c. forcing over $M$,
because one can simply add a Cohen real, but it is not true in
$M$ by the very definition of $A$. More generally, if $A$ is any
set in $M$ having hereditary size at least $2^\omega$, then the
assertion, ``{\it $2^\omega$ is greater than the hereditary size
of $A$},'' is forceably necessary by c.c.c. forcing, since one
may make $2^\omega$ as large as desired, but again, it is not
true in $M$.

This argument suggests that we must at least restrict our
parameter space to $H(2^\omega)$, those parameters of hereditary
size less than $2^\omega$. And the following theorem shows that
restricting this far is enough.

\Theorem. The following theories are equiconsistent.
\begin{enumerate}
\item $\ZFC+\MPtilde_{ccc}$, allowing arbitrary parameters
from $H(2^\omega)$, plus $2^\omega$ is weakly inaccessible.
\item $\ZFC+\MPtilde_{ccc}$, allowing ordinal parameters below
$2^\omega$, plus $2^\omega$ is regular.
\item $\ZFC+{V_\delta\elesub V}+\delta$ is inaccessible.
\item $\ZFC+\MPtilde$.
\end{enumerate}\label{MPccc}

\Proof: Clearly, the first theory directly implies the second
theory. Furthermore, I have already shown that the third and
fourth theories are equiconsistent.

\SubLemma. Every model of theory 2 contains an inner model of
theory 3.

This Lemma is a consequence of the following Lemma:

\SubLemma. If $\MPtilde_{ccc}$ holds (in the language with
ordinal parameters below $2^\omega$) and $\delta=2^\omega$ is
regular, then $\delta$ is inaccessible in $L$ and $L_\delta\elesub
L$.\label{innerccc}

\Proof: Suppose that the Maximality Principle $\MPtilde_{ccc}$
holds for c.c.c.~forcing in the language with ordinal parameters
below $\delta=2^\omega$, and that $\delta$ is regular. Note that
for any $\gamma<\delta$, the assertion
``$2^\omega>\gamma^\plus$'' is forceably necessary, since by
c.c.c.~forcing, one can pump up $2^\omega$ as large as desired.
Thus, this assertion must be true, and so $\delta$ is weakly
inaccessible and hence inaccessible in $L$.

To prove $L_\delta\elesub L$, I will simply verify the
Tarski-Vaught criterion. Suppose that $L\satisfies\exists
x\,\psi(x,y)$ for some $y\in L_\delta$. Let $\beta<\delta$ be
large enough so that $y\in L_\beta$. Now, by the Replacement
Axiom, let $\alpha$ be least so that for every $y'\in L_\beta$
with an $x$ satisfying $\psi(x,y)$, there is such an $x$ in
$L_\alpha$. This property defines $\alpha$ in any forcing
extension, and so the assertion $2^\omega>\alpha$ is forceably
necessary. Thus, $\alpha<\delta$, and so the witness $x$ for $y$
can be found in $L_\delta$, as desired.\QED

\SubLemma. Any model of theory 3 has a forcing extension that is
a model of theory 1.\label{outerccc}

\Proof: This argument follows Lemma \ref{outer} above. Suppose
that $V_\delta\elesub V$ and $\delta$ is inaccessible. Enumerate
$\<\varphi_\alpha\st\alpha<\delta>$ all sentences $\varphi$ (with
unbounded repetition) in the language of set theory with (names
for) parameters in $V_\delta$ coming from a forcing extension by
forcing of size less than $\delta$. Define a finite support
$\delta$-iteration of c.c.c.~forcing, so that at each stage
$\gamma<\delta$, the forcing $\Q_\alpha$ forces $\varphi_\alpha$
to be necessary over $V_\delta^{\P_\alpha}$, if possible (i.e.
first, if this makes sense, if the parameters appearing in
$\varphi_\alpha$ are $\P_\alpha$-names, and second, if there is
such a forcing in $V_\delta^{\P_\alpha}$). Let $G\of\P_\delta$ be
$V$-generic for $\P_\delta$, and consider the model $V[G]$. Note
that periodically during the iteration, the posets will ensure
that $2^\omega$ is made arbitrarily large below $\delta$, and so
$2^\omega\geq\delta$ in $V[G]$. Conversely, since the forcing has
size $\delta$, it is easy to see that $V[G]\satisfies
2^\omega=\delta$. Further, since the iteration is c.c.c.,
$\delta$ remains a regular limit cardinal. Finally, if $\varphi$
is a forceably necessary assertion over $V[G]$ with parameters in
$H(2^\omega)=H(\delta)$, then $\varphi=\varphi_\alpha$ for some
$\alpha$ such that the parameters are $\P_\alpha$-names. Since
$\varphi$ was forceably necessary in $V[G]$ it must have been
forceably necessary in $V[G_\alpha]$ and hence, since
$V_\delta[G_\alpha]\elesub V[G_\alpha]$, it was forceably
necessary in $V_\delta[G_\alpha]$. Thus, it was forced to be
necessary in $V_\delta[G_{\alpha+1}]$, and so it was necessary in
$V[G_{\alpha+1}]$. Thus, it is true in $V[G]$, as desired. In
summary, $V[G]\satisfies\MPtilde_{ccc}+2^\omega$ is weakly
inaccessible.\QED

This completes the proof of the theorem.\QED

A nearly identical argument works in the case of proper forcing:

\Theorem. The following theories are equiconsistent.
\begin{enumerate}
\item $\ZFC+\MPtilde_{\rm proper}$, using parameters in
$H(2^\omega)$, plus $2^\omega$ is weakly inaccessible.
\item $\ZFC+\MPtilde_{ccc}$, using parameters in $H(2^\omega)$, plus $2^\omega$ is weakly inaccessible.
\item $\ZFC+\MPtilde$.
\end{enumerate}
\label{MPproper}

\Proof: Using the techniques of \ref{innerccc}, one can show that
if the first theory holds, then $\delta=2^\omega$ is inaccessible
in $L$ and $L_\delta\elesub L$. And this has already been proved
equiconsistent with the second and third theories.

Conversely, the second and third theories are equiconsistent with
$\ZFC+V_\delta\elesub V$ for an inaccessible cardinal $\delta$.
Given this, one enumerates the formulas $\<\varphi_\alpha\st
\alpha<\delta>$ using (names for) parameters in $V_\delta$
appearing in forcing extensions of size less than $\delta$, and
then performs a countable support $\delta$-iteration $\P_\delta$
of proper forcing, which at stage $\alpha<\delta$ forces the
necessity of $\varphi_\alpha$ over $V_\delta^{\P_\alpha}$ by
proper forcing, provided that this make sense (that is, the names
in $\varphi_\alpha$ are $\P_\alpha$-names) and that this is
possible (that is, there is a proper forcing poset accomplishing
this). The resulting forcing extension $V^\P$ will be a model of
$\ZFC+\MPtilde_{\rm proper}$, just as in Lemma \ref{outerccc}.\QED

And by using revised countable support in the iteration, one can
deduce a similar fact for $\MPtilde_{semi-proper}$, using
parameters in $H(2^\omega)$.

A natural question remains open here, whether the assertion that
$2^\omega$ is weakly inaccessible can be removed from Theorems
\ref{MPccc} and \ref{MPproper}. The arguments I gave show that
$\MPtilde_{ccc}$ implies that $2^\omega$ is a limit cardinal. And
it is easy to see that it must be a limit cardinal of very high
rank, since with c.c.c.~forcing one can push $2^\omega$ beyond
$\aleph_{\omega_1}$, $\aleph_{\aleph_{\omega_1}}$ or any other
definable cardinal. Does $\MPtilde_{ccc}$ imply outright that
$2^\omega$ is regular? If so, then $\MPtilde_{ccc}$ is
equiconsistent with $\MPtilde$. If not, then is $\MPtilde_{ccc}$
actually weaker than $\MPtilde$ in consistency strength?

Another interesting point is that while the arguments here are
based on that of Lemma \ref{outer}, in that earlier lemma we could
actually say much more about the nature of the $\delta$-iteration
$\P$, namely, that the forcing $\P$ there was necessarily the
L\'evy collapse making $\delta$ the $\omega_1$ of the extension.
In the previous argument and in Lemma \ref{outerccc}, this is
certainly not the case. But is there, nevertheless, an
alternative, simple characterization of the iterations $\P$? If
so, this would be very interesting to know.

There are truly a plethora of open problems remaining here. Apart
from the big open question---the relative consistency of
$\necessary\MPtilde$---one can analyze the various modified
maximality principles, including $\MP_{\rm card}$,
$\MP_{\omega_1}$, $\MP_{\rm stat-pres}$, and so on, along with
their boldface and necessary forms ($\MPtilde_{\rm card}$,
$\necessary\MPtilde_{\rm card}$, etc.). Other interesting
principles are obtained by restricting to forcing extensions that
preserve a given statement or theory, such as $\CH$. For example,
$\MP_{\hbox{\scriptsize\sc ch}}$ is the scheme asserting that
every statement that is forceably necessary by forcing preserving
\CH\ is already true. One can also consider forcing notions not
adding reals, or preserving a given supercompact cardinal, and so
on. The possibilities seem endless.

\end{document}